\documentclass[11pt]{article}

\usepackage{authblk}
\usepackage{amsfonts}
\usepackage{amsthm}
\usepackage{amsmath}
\usepackage{graphicx}
\usepackage{comment}
\usepackage[usenames,dvipsnames]{xcolor}
\usepackage{enumerate}
\usepackage{extarrows}
\usepackage{hyperref}
\usepackage[round, authoryear]{natbib}

\usepackage{stackrel}
\usepackage{centernot}
\usepackage{caption}
\usepackage{subcaption}
\usepackage{fancyhdr}
\let\quoteOLD\quote
\def\quote{\quoteOLD\small}

\definecolor{labelkey}{cmyk}{0,0.8,1,0.5}
\definecolor{refkey}{cmyk}{0,0.8,1,0.5}

\makeatletter
\def\th@newremark{\th@remark\thm@headfont{\bfseries}}
\makeatletter

\def\boxit#1{\vbox{\hrule\hbox{\vrule\kern6pt
          \vbox{\kern6pt#1\kern6pt}\kern6pt\vrule}\hrule}}

\newtheorem{theorem}{Theorem}

\newtheorem{proposition}{Proposition}
\newtheorem{definition}{Definition}

\newtheorem{remark}{Remark}

\numberwithin{equation}{section}
\numberwithin{theorem}{section}
\numberwithin{corollary}{section}
\numberwithin{proposition}{section}
\numberwithin{lemma}{section}
\numberwithin{definition}{section}
\numberwithin{remark}{section}

\newcommand{\lambar}{\overline{\Lambda}}
\newcommand{\laminv}{\overline{\Lambda}^{\leftarrow}}

\newcommand{\eqdr}{\stackrel{\mathrm{D}}{=}}

\newcommand{\eqd}{\stackrel{{\mathrm D}}=}

\newcommand{\R}{\Bbb{R}}

\newcommand{\N}{\Bbb{N}}

\newcommand{\rmd}{{\rm d}}

\newcommand{\halmos}{\quad\hfill\mbox{$\Box$}}

\newcommand{\FFF}{{\cal F}}

\newcommand{\BB}{\mathbb{B}^{(r)}}
\newcommand{\BBi}{\mathbb{B}^{(r+i)}}

\newcommand{\BBy}{\mathbb{B}^{(r+1)}}
\newcommand{\BBz}{\mathbb{B}^{(r+2)}}
\newcommand{\BBx}{\mathbb{B}^{(r,x-)}}
\newcommand{\BX}{\mathbb{B}^{(r,x_1-)}}

\newcommand{\MM}{\mathbb{M}}
\newcommand{\NN}{\mathbb{N}}

\newcommand{\wwt}{\widetilde}
\newcommand{\wt}{\tilde}

\newcommand{\JJJJ}{\mathfrak{J}}

\newcommand{\be}{\begin{equation}}
\newcommand{\ee}{\end{equation}}
\newcommand{\bea}{\begin{eqnarray}}
\newcommand{\eea}{\end{eqnarray}}
\newcommand{\bean}{\begin{eqnarray*}}
\newcommand{\eean}{\end{eqnarray*}}
\newcommand{\ben}{\begin{equation*}}
\newcommand{\een}{\end{equation*}}
\newcommand{\ba}{\begin{aligned}}
\newcommand{\ea}{\end{aligned}}

\def\nexto{\kern -0.54em}

\newcommand{\PP}{\textbf{\rm P}}
\newcommand{\EE}{\textbf{\rm E}}
\newcommand{\PD}{\textbf{\rm PD}}
\newcommand{\BN}{\mathcal{BN}}

\title{Generalised Poisson-Dirichlet  Distributions and the Negative Binomial Point Process}

\author{ Yuguang F. Ipsen\thanks{ Y. F. Ipsen thanks the support from the ARC Center of Excellence for Mathematical and Statistical Frontiers. Corresponding author: yuguang.fan@anu.edu.au} }

\author{ Ross A. Maller\thanks{Research partially supported by ARC Grant DP1092502; Email: Ross.Maller@anu.edu.au}}
\affil{ Research School of Finance,  Actuarial Studies \&  Statistics\\
 Australian National University, Australia.}

\begin{document}

\maketitle

\begin{abstract}
When  $S=(S_t)_{t\ge 0}$ is an $\alpha$-stable subordinator, the sequence of ordered jumps of $S$, up till time $1$, omitting the $r$ largest of them,  and taken as  proportions  of  their sum $^{(r)}S_t$,   defines a 2-parameter distribution on the infinite dimensional  simplex, $\nabla_{\infty}$, which we call the $\PD_\alpha^{(r)}$ distribution. 
When $r=0$ it reduces to the $\PD_\alpha$ distribution introduced by Kingman in 1975. 
We observe a serendipitous connection between  $\PD_\alpha^{(r)}$ and the negative binomial point process of
Gregoire (1984), which we 
 exploit to analyse in detail a size-biased version of  $\PD_\alpha^{(r)}$.
As a consequence we derive a  stick-breaking representation for
 the process and a useful form for its distribution.
  This program  produces a large new class of distributions available for  a variety of  modelling purposes.

\end{abstract}



\section{Introduction}\label{intro}
Developments related to the Poisson-Dirichlet distribution and its generalisations have had an enormous impact in recent times, stimulating as well as synthesising a host of theoretical results connected in  particular to the excursion theory of stochastic processes and to random partitions, and opening up a wealth of applications areas, especially for example in Bayesian statistics and population genetics. 
 We refer to Bertoin  \citeyearpar{ber96} and Feng \citeyearpar{Feng2010}
for up-to-date accounts of various aspects.

To motivate the ideas that concern us here, start with a stable 
subordinator  $(S_t)_{t\ge 0}$ of index $\alpha\in(0,1)$ 
on $\R^+$ having   jump process $(\Delta S_t:=S_t-S_{t-})_{t>0}$, and order the jumps up till 1 as $\Delta S_1^{(1)}\ge \Delta S_1^{(2)}\ge\cdots$.
The random sequence $(\Delta S_1^{(i)}/S_1)_{i\ge 1}$ specifies a distribution on the infinite dimensional simplex $\nabla_{\infty}$ which we will refer to as a $\PD_\alpha$ distribution. It was introduced by \cite{kingman1975} and subsequently gave rise to a large body of research. Of special interest to us are papers by  
\cite{PPY1992} and \cite{PY1992, PY1997} (hereafter, referred to as PPY \citeyearpar{PPY1992} and PY \citeyearpar{PY1992, PY1997}). They contain in particular formulae for  the distribution of the {\it size-biased} vector associated with
 $\PD_\alpha$.

The $\PD_\alpha$ distribution arises by considering the ordered jumps $(\Delta S_1^{(i)})_{i\ge 1}$ and their relation to the sum, $S_1$.
As a natural generalisation,  delete the $r$ largest jumps ($r\in\NN$)\footnote{Throughout, let $\N:=\{1,2,\ldots\}$ and $\N_0:=\{0,1,2,\ldots\}$.} up till time $1$ and consider the distribution of the remaining jumps $(\Delta S_1^{(i)})_{i\ge r+1}$ taken as proportions of their sum, $^{(r)}S_1$, the latter being $S_1$ with the $r$ largest terms removed. Again we obtain
 a distribution on $\nabla_{\infty}$, now with an extra parameter, $r$.
 When $r=0$ (no trimming) this is a $\PD_\alpha$ distribution, while for $r=1,2,\ldots$, it defines a 2-parameter distribution on  $\nabla_{\infty}$ which we call the $\PD_\alpha^{(r)}$ distribution.

Laplace transforms of the stable ratios take a reasonably explicit form, and 
reveal a close  connection with the {\it negative binomial point process} of \cite{Gregoire1984}.
This suggests  some rewarding new lines of enquiry, and we proceed to define a size-biased version of $\PD_\alpha^{(r)}$ and use the point process representation to derive a corresponding stick-breaking representation.
  
  This program  produces a large new class of distributions 
  available for  a variety of  modelling purposes.
  We illustrate its relevance by reference to two applied situations, occurring in finance and linguistics which are analysed in papers by \cite{sos15} and \cite{ggj11}.
We discuss these further in Section \ref{dis}.

\section{Jumps of a Normalised Stable Subordinator}\label{GPD}
We consider a driftless stable subordinator, that is, a  real valued L\'evy process $(S_t)_{t\ge 0}$, with $S_0\equiv 0$, 
on a filtered probability space $(\Omega, ({\cal F}_t)_{t\ge 0}, \PP)$,  with canonical triplet $(\gamma,0,\Lambda)$;  thus, 
 having L\'evy measure 
\be\label{LM2}
\Lambda(\rmd x)=c\alpha x^{-\alpha-1}\rmd x{\bf 1}_{\{x>0\}}, \quad \text{for some }c>0 \ \text{and }  0<\alpha<1,
\ee
with tail measure
\ben
\lambar(x)=c x^{-\alpha}, \quad \quad x > 0,
\een
and Laplace transform 
 \be\label{psdef}
 Ee^{-\lambda S_t}= e^{-t\Psi(\lambda)}, \quad \text{where} \quad 
\Psi(\lambda) =
\int_{(0,1)} \big(1-e^{-\lambda x}\big)  \Lambda(\rmd x),\
\lambda>0.
\ee
 Write $(\Delta S_t:=S_t-S_{t-})_{t>0}$, with   $\Delta S_0=0$, for the jump process of $S$, and
 $\Delta S_t^{(1)} \ge \Delta S_t^{(2)}\ge \cdots$
 for the ordered jumps at time $t>0$.
 Since $\Lambda\{(0,\infty)\}=\infty$ there are infinitely many jumps, a.s. (almost surely), in any finite time interval, and since   $\Lambda$ is diffuse,
 the  ordered jumps $(\Delta S_t^{(i)})_{i=1,2,\ldots}$ are uniquely defined, a.s. 
 Our objective 
 is to study the 
 ``trimmed" process, by which we mean $S_t$ minus its large 
 jumps, at a given time $t$:
\begin{equation}\label{trims}
{}^{(r)}S_t:= S_t-\sum_{i=1}^r {\Delta S}_t^{(i)}, \ r\in\N, \ t>0
\end{equation}
(and we set $^{(0)}S_t\equiv S_t$).

In the next subsection  we define the 
$\PD_\alpha^{(r)}$ distribution
as that of the sequence of ratios $\big(\Delta S_1^{(i)}/\,{}^{(r)}S_1\big)_{i\ge r+1}$.

\subsection{Generalised Poisson-Dirichlet  Distributions}\label{ssGPD}
Fix $r\in\N_0$ and define
\ben\label{defVS}
V_n^{(r)}:=  \frac{\Delta S^{(r+n)}_1}{{}^{(r)}S_1},\ n\in\N.
\een
%
Since, for $m>n$,
 \ben
\sum_{n<j\le m} V_j^{(r)}\le 
\frac{{}^{(r+n)}S_1}{^{(r)}S_1} \to 0,\ {\rm a.s.,\ as}\ n\to \infty,
\een
 the series  $\sum_{n} V_n^{(r)}$ converges a.s. for each $r\in\N$, and clearly,  $\sum_{n\ge 1} V_n^{(r)}=1$.
 Consequently, the distribution of  $(V_n^{(r)})_{n\in\N}$  when $r\in\N$ defines a  new family of distributions on $\nabla_{\infty}$
derived from  the subordinator $S$.

\begin{definition}\label{deff}
Let $(S_t, 0\le t \le 1)$ be a driftless stable subordinator with index $\alpha\in(0,1)$ and take $r\in\N_0$. Then the distribution of the sequence 
\be\label{sub_alpha}
\big(V_n^{(r)}\big)_{n\in\N} = \big( V_1^{(r)},\,  V_2^{(r)}, \ldots \big)=  \Big( \frac{\Delta S_1^{(r+1)}}{{}^{(r)}S_1},\, \frac{\Delta S_1^{(r+2)}}{{}^{(r)}S_1}, \ldots \Big)
\ee 
we call  a $\PD_\alpha^{(r)}$ distribution. When $r =0$, $\PD_\alpha$ is recovered.
\end{definition}

\begin{remark}\label{def3e}
{\rm 
$\PD_\alpha^{(r)}$ is obtained from the deletion of the $r$ largest jumps of $S_1$, followed by renormalisation, 
rather than  from the deletion of the first $r$ size-biased picks from $(\Delta S_1^{(i)})$, as considered in \cite{Pitman2003} and PY \citeyearpar[Prop. 34, 35]{PY1997}. 
This results in a different dependence structure in $\PD_\alpha^{(r)}$ for the stick-breaking representation than in the PY situations
(cf.   Theorem \ref{stick} in the next section).

Results related to those of PY \citeyearpar{PY1997} concerning deletion of excursion intervals of certain Bessel bridges are in PPY \citeyearpar[Sect.3]{PY1992}; see also James \citeyearpar{James2013, James2015}.
}
\end{remark}

\begin{remark}\label{def3}
\noindent{\rm 
Similar to \eqref{sub_alpha}, any distribution on $\nabla_\infty$ with a subordinator representation can be generalised by removing the $r$ largest jumps up till time $t>0$ from the subordinator.
For example: \ (i) \ 
the usual Poisson-Dirichlet distribution, denoted as  
$\PD(0, \theta)$ in PY~\citeyearpar{PY1997}, 
can be generalised by trimming a Gamma subordinator up till time $\theta>0$;

(ii) \ the two parameter 
$\PD(\alpha, \theta)$  distribution 
in PY~\citeyearpar{PY1997} can be extended by trimming a generalised Gamma subordinator up till a random time mixed with a Gamma$(\theta/\alpha,1)$ distribution (see PY~\citeyearpar[Prop. 21]{PY1997}).
}
\end{remark}

We do not pursue these generalisations here, going on instead to explore a connection with the negative binomial process (in the next subsection), but we conclude this subsection with a formula for 
the Laplace transform of the ratio  $^{(r)}S_1/\Delta S_1^{(r)}$:
\be\label{31L}
\EE\big(e^{-\lambda{}^{(r)}S_1/\Delta S_1^{(r)}}\big)
=\EE\big(e^{-\lambda W_{\Gamma_{r}}}\big)
=\frac{1}{\left(1+\Psi(\lambda)\right)^{r}}, \ r\in\N,
\ee
where $W=(W_v)_{v\ge 0}$ is a  driftless subordinator with measure 
$\Lambda(\rmd x){\bf 1}_{(0,1)}$, and $\Gamma_r$ denotes a Gamma$(r,1)$ random variable independent of $W$. Recall that $\Psi$ is defined in \eqref{psdef}.
Formula \eqref{31L} appears in PY~\citeyearpar[Prop. 11]{PY1997}, as well as in \cite{KeveiMason2014}.

Later we will need the density function of a Gamma  random variable  with parameter $r$:
\ben
\PP(\Gamma_r\in\rmd x)= \frac{x^{r-1} e^{-x}\rmd x}{\Gamma(r)} {\bf 1}_{\{x>0\}},
\een
and the density of a Beta  random variable $B_{a, b}$ with parameters $a, b>0$:
\be\label{betadef}
 f_B(x)= \frac{\Gamma(a+b)}{\Gamma(a)\Gamma(b)} x^{a-1}(1-x)^{b-1}{\bf 1}_{\{0<x<1\}},
 \ee
where $\Gamma(r)=\int_0^\infty x^{r-1}e^{-x}\rmd x$, $r>0$,  is the Gamma function.


\subsection{The Negative Binomial Point Process and the Distribution of
$\PD^{(r)}_\alpha$}\label{ssNB}
In this subsection we connect the previous results 
with  the negative binomial point process introduced by \cite{Gregoire1984}.
As previously, $S$ is the driftless stable subordinator with index $\alpha\in(0,1)$. 
It turns out that to generate the Laplace transform in \eqref{31L}, we have to construct a point process from ratios of stable jumps, rather than from the jumps themselves. Thus, for $r\in\N$, and with $\delta_x$ denoting a point mass at $x\in\R$, 
define a random point measure on the Borel sets of $(0,1)$ by 
\be\label{defB}
 \BB = \sum_{i \ge 1} \delta_{J_r(i)}, \quad \text{where} \quad 
 J_r(i) = \frac{\Delta S_1^{(r+ i)}}{\Delta S_1^{(r)}},\ i = 1, 2, \ldots.
\ee


Let $(\mathbb{M},{\cal M})$
be the space of all point measures\footnote{We generally follow the exposition in Resnick \citeyearpar[Chap. 3]{res87} for the following setup.} on $(0,1)$ with its usual Borel $\sigma$-algebra
and let $\FFF^+$ be the set of nonnegative measurable functions on $(0,1)$.
A random measure $\xi$ on $(\mathbb{M}, {\cal M})$ has Laplace functional
defined as
\ben 
 \Phi(f) = \EE (e^{-\xi(f) }) = \int_{M\in\mathbb{M}} 
e^{-\int_{(0,1)}f(x)M(\rmd x)} \PP(\xi\in \rmd M), \ f \in \FFF^+.
\een
 Given a measure $\Pi$ on $(0,\infty)$,  locally finite at infinity, \cite{Gregoire1984} defines the point process
 $\BN(r, \Pi)$ on $(\mathbb{M}, {\cal M})$
in terms of its Laplace  functional as
\be\label{NBp}
 \Phi(f) =\Big(1 + \int_{(0,\infty)}(1-e^{-f(x)})\Pi(\rmd x) \Big)^{-r}, \ f \in \FFF^+.
\ee
Recall \eqref{LM2} and let $\wt \Lambda(\rmd x) : = \alpha x^{-\alpha-1}\rmd x {\bf 1}_{\{0<x<1\}}$ be the normalised $\Lambda(\rmd x)$ restricted to $(0,1)$.
For $r\in\N$, denote the law of $\BN(r, \wt \Lambda)$ by $\PP_r$, so that 
\[
\PP_r(\rmd M)=\PP\big(\BN(r, \wt\Lambda)\in \rmd M\big),  \quad M\in \MM. 
\]
Let the family of Palm distributions of $\PP_r$ be $\PP_r^{(x)}$, $x \in (0,1)$.
The correspondence between \eqref{NBp} and the righthand side of \eqref{31L} suggests the following proposition. 

\begin{proposition}\label{BN0}
Let $\BB$ be defined as in \eqref{defB}. Then

{\rm (i)}
$\BB$ is a negative binomial point process with distribution $\PP_r$ such that 

{\rm (ii)}
$ \EE(\BB(A)) = r \wt \Lambda(A)$  for any Borel set $A \subset (0,1)$.

{\rm (iii)}
The Laplace functional of the probability measure $\PP_r^{(x)}$ on $(\mathbb{M}, {\cal M})$ satisfies
\be\label{palmBN}
\Phi_{\PP_r^{(x)}}(f)  = \Phi_{\delta_{x}}(f) \Phi_{\PP_{r+1}}(f), \quad  f \in \FFF^+.
\ee
\end{proposition}

\begin{remark}[Interpretation of  $\PP_r^{(x)}$]\label{rem3c}
{\rm 
We can think of the Palm distribution $\PP_r^{(x)}$ as the conditional distribution of 
 $\BN(r, \wt \Lambda)$ given  $\BN(r, \wt \Lambda)(\{x\})>0$. 
From \eqref{palmBN}, we can interpret 
$\PP_r^{(x)}$ in the following way. Let $\xi $ be distributed as $\BN(r+1, \wt \Lambda)$. Then $\PP_r^{(x)}$ is the distribution of $\xi+ \delta_{x}$.
}
\end{remark}

\medskip\noindent {\bf Proof of Proposition \ref{BN0}:}\
(i):\
Conditional on $\{\Delta S_1^{(r)} = v\}$, $v>0$, the truncated point process $\{\Delta S_1^{(r+j)},\, j\in\N\}$ is a Poisson point process with intensity measure $\Lambda(\rmd x)1\{x < v\}$. 
The following distributional equivalence can be deduced, for example, from Lemma 1.1 of \cite{bfm2016}:
\be\label{ch}
\big(\Delta S_t^{(i)}\big)_{1\le i\le r}
\eqdr \big(\laminv(\Gamma_i/t)\big)_{1\le i\le r},\ t>0,\ r\in\N,
\ee
where the $\Gamma_i$ are Gamma$(i,1)$ random variables and $\laminv(x) = c^{1/\alpha}x^{-1/\alpha}$ is the inverse function of $\lambar$.
Using \eqref{ch} we can write the Laplace functional of $\BB$ as 
\begin{align}\label{BN2}
\EE(e^{-\BB(f)}) &= 
\int_{v>0}\exp\bigg( -\int_{0}^{v} (1-e^{-f(x/v)}) \Lambda(\rmd x)\bigg) \PP\big(\laminv(\Gamma_r) \in \rmd v\big) \nonumber \\
&=
\int_{v>0}\exp\bigg( -\int_{0}^{1} (1-e^{-f(x)}) \Lambda(v \rmd x)\bigg) \PP\big(\laminv(\Gamma_r) \in \rmd v\big) \nonumber \\ 
 &=
 \int_{v>0} \exp\bigg( - cv^{-\alpha}\int_{0}^{1} (1-e^{-f(x)}) \alpha x^{-1-\alpha} \rmd x\bigg) \PP\big(\Gamma_r \in \rmd (c v^{-\alpha})\big),
\end{align}
for each $f\in \FFF^+$.
By change of variable in \eqref{BN2} with $y = c v^{-\alpha}$, we have
\begin{align}\label{BN3}
\EE(e^{-\BB(f)})&=\int_{0}^\infty  \exp\bigg( -y\int_{0}^{1} (1-e^{-f(x)}) \wt \Lambda(\rmd x)\bigg) \PP(\Gamma_r \in \rmd y)\nonumber \\
		        &=\int_{0}^\infty  \exp\bigg( -y\int_{0}^{1} (1-e^{-f(x)}) \wt \Lambda(\rmd x)\bigg)  \frac{y^{r-1}e^{-y}}{\Gamma(r)}\rmd y \nonumber \\ 
		        &= \bigg( 1 +  \int_{0}^{1} (1-e^{-f(x)}) \wt \Lambda(\rmd x) \bigg)^{-r}.
\end{align}
Comparing \eqref{BN3} with \eqref{NBp} proves Part (i).
Parts {\rm (ii)} and (iii) follow from Part {\rm (i)} by Propositions 3.3 and 4.3 in \cite{Gregoire1984}. 
\halmos

\begin{remark}\label{def3a}
{\rm 
(i)\ The sum of the points in $\BB$ is $^{(r)}S_1/\Delta S_1^{(r)}$, hence the connection with \eqref{31L}. 

(ii)\ 
A variety of formulae relating to the Poisson-Dirichlet distributions have been derived over the years, including an iterative formula for the joint density of the first $n$ terms of $\PD_\alpha$
(Perman \citeyearpar[Thm. 2]{perman1993}).
Such formulae, while  
explicit, are ``rather intractable"  
(PY \citeyearpar[p.329]{PY1992}), and simpler structures can be revealed for the
corresponding  size-biased permutation;
see, e.g., PPY \citeyearpar[Thm. 1.2]{PPY1992} (attributed to \cite{perman_thesis}), which allows for a ``stick-breaking" representation of $\PD_\alpha$ in terms of independent Beta rvs. 
See also \cite{PitmanTran2015} for the setting of a finite sequence of i.i.d. random variables.
This motivates us to consider the size-biased permutation of $\PD_\alpha^{(r)}$ and to investigate a stick-breaking-like representation in the $r$-trimmed case through
the random point measure $\BB$.

\cite{Pitman1995_SSB} proved that the $\PD(\alpha,\theta)$ of PY \citeyearpar{PY1997} is the largest class of distributions with 
a stick-breaking representation in terms of independent beta rvs; inevitably, then, our enlarged class $\PD_\alpha^{(r)}$ 
requires a dependent  stick-breaking representation. (\cite{James2013} derives another class of distributions, PG$(\alpha, \zeta)$, by mixing generalised Gamma subordinators, which also has a dependent stick-breaking representation.) The dependence structure will become clear in the main  Theorem \ref{stick} of this section  which gives a formula for the density of the size-biased version of the sequence $(V_n^{(r)})$ in \eqref{sub_alpha}. 
The remaining calculations in this section lead up to  Theorem \ref{stick}.
}
\end{remark}

Henceforth fix $r \in \N$.
Write 
\be\label{Jd}
\JJJJ_r : = \{ J_r(1), J_r(2), J_r(3), \ldots\}
\ee
for the points of $\BB$, 
with sum 
\be\label{rT}  
{}^{(r)}T:= \sum_{i \ge 1} J_r(i) = \frac{{}^{(r)}S_1}{\Delta S_1^{(r)}}.
\ee
Define the \emph{size-biased random permutation} of $\JJJJ_r$, denoted by  $ ( \wt J_1, \wt J_2, \wt J_3, \ldots )$, in the following way.  Conditional on $\JJJJ_r$, the first term $\wt J_1$ takes values among the members of  $\JJJJ_r$ with probabilities
\ben\label{sb1}
 \PP\big( \wt J_1 =  J_r(i) \, \big | \, \JJJJ_r \big) 
 =\frac{J_r(i)}{\sum_{\ell\ge 1}J_r(\ell)}
 = \frac{\Delta S_1^{(r+i)}}{ ^{(r)}S_1}, \ i=1,2,\ldots
\een
 Conditional on $\JJJJ_r$ and $ \wt J_1, \ldots,  \wt J_{n}$, for each $n = 1,2, \ldots$, the $(n+1)^{st}$ term   $\wt J_{n+1}$ takes values among 
 $\{ J_r(j), j=1,2,\ldots \, ; \,   J_r(j) \neq  \wt J_l, 
 1\le l \le n\}$, with probabilities
\ben\label{sb2}
\PP( \wt J_{n+1}=  J_r(j) \big |  \wt J_1, \ldots,  \wt J_{n}, \JJJJ_r) = \frac{\Delta S_1^{(r+j)}1\{ \Delta S_1^{(r+j)} \neq  \wt J_l  \cdot \Delta S_1^{(r)},1\le l \le n\}}
{^{(r)}S_r- \Delta S_1^{(r)} \cdot (\sum_{l=1}^{n} \wt J_{l}) }.
\een
Then the sums of the  remaining points in the point process,  
after removing points by size-biased sampling, are
\be\label{sb_T}
 {}^{(r)}T_1 := {}^{(r)}T - \wt J_1, \quad \text{and for each } n> 1, \quad {}^{(r)}T_n := \, {}^{(r)}T_{n-1} - \wt J_n.
\ee
The successive {\it residual fractions} are
\be\label{U1}
 ^{(r)}U_1 := \frac{^{(r)}T_1}{{}^{(r)}T}
 =1-\frac{\wt J_1}{^{(r)}T},
 \ee
 and for each $n >1$, 
 \be\label{Un}
^{(r)}U_n := \frac{^{(r)}T_n}{^{(r)}T_{n-1}}
 =1-\frac{\wt J_n}{^{(r)}T_{n-1}}.
\ee

For $M$ a point measure in $\mathbb{M}$, let
$T(M) = \sum_{x \in M} x$ be the sum of the magnitudes of the points in $M$.
For each $r \in\N$ let the density of $T(\mathbb{B}^{(r)})$ with $\mathbb{B}^{(r)}$ distributed as $\BN(r, \wt \Lambda)$ be
\be\label{gd}
g_{r}(t) := \PP\big( T(\BB) \in \rmd t\big)/\rmd t = \PP\big({}^{(r)}T \in \rmd t\big) /\rmd t,\ t>0.
\ee
By \eqref{rT}, ${}^{(r)}T= {}^{(r)}S_1/\Delta S_1^{(r)}$, so by \eqref{31L}, $g_{r}$ satisfies
\be\label{gdL}
\int_0^\infty e^{-\lambda x} g_{r}(x) \rmd x = 
\Big(1 + \int_{0}^{1} (1- e^{-\lambda x}) \wt \Lambda(\rmd x)  \Big)^{-r},
\ee
 for $\lambda > 0$, $r\in\N$. 
Alternatively expressed, $g_{r}$ is the density of $\wwt W_{\Gamma_{r}}$, where $(\wwt W_v)_{v\ge 0}$ is the driftless subordinator with L\'evy measure  $\wt\Lambda(\rmd x)$.

The next lemma derives important properties of $\BB$. 
It will be apparent that our proofs owe much to the methods of PPY \citeyearpar{PPY1992}, PY \citeyearpar{PY1992} and Fitzsimmons, Pitman and Yor \citeyearpar[Sect. 5]{Fitzsimmons1993}.
See  Lemma 2.2 of PPY \citeyearpar{PPY1992}.
 In a remark at the end of this section  we discuss briefly the differences as well as similarities between our
  approaches.

\begin{proposition}\label{b11}
For $r\in\N$ let $\BB$ be a negative binomial point process. 
Let $\BB_1 = \BB - \wt J_1$ be the remaining process after removing the first size-biased pick. 

{\rm (i)}\ Then for $0<x<1$, $M\in\mathbb{M}$, 
\be\label{b1}
\PP(\wt J_1 \in \rmd x, \, \BB_1 \in \rmd M) 
= \frac{x}{T(M) + x} r \wt \Lambda(\rmd x ) 
\PP\big(\BBy\in \rmd M\big).
\ee

{\rm (ii)}\  For $0<x<1$, $M\in\mathbb{M}$, $t>0$,
we have 
\begin{align}\label{b1T}
&\PP\big(\wt J_1 \in \rmd x, \, \BB_1\in\rmd M,\, {}^{(r)}T_1 \in \rmd t\big) \nonumber \\
&= \frac{x}{t + x} r \wt \Lambda(\rmd x ) 
\PP\big(\BBy\in \rmd M, \, T(\BBy) \in \rmd t \big).
\end{align}

{\rm (iii)}\ 
For $0<x<1$, $t>0$, we have
\be\label{b1T0}
\PP(\wt J_1 \in \rmd x, \,   {}^{(r)}T_1  \in \rmd t  ) 
= \frac{x}{t + x} r \wt \Lambda(\rmd x ) \PP \big( {}^{(r+1)}T \in \rmd t \big).
\ee
\end{proposition}

\medskip\noindent {\bf Proof of Proposition \ref{b11}:}\
(i)\ 
The definition of the size-biased picks implies
\be\label{b1a}
\PP(\wt J_1 \in \rmd x \, | \,  \BB=M)=  \frac{x}{T(M)}M(\rmd x),
\ 0<x<1,\ M \in \mathbb{M}\setminus\{\phi\}.
\ee
Recall that $\PP_r(\rmd M) = \PP(\BB \in \rmd M)$.
We use  the following property of \emph{Palm distributions} (see for instance Daley and Vere-Jones \citeyearpar[Sect. 12.1]{DaleyVereJones1998}):
\be\label{palm0}
r \wt \Lambda(\rmd x) \PP^{(x)}_r(\rmd M) = M(\rmd x) \PP_r(\rmd M) = M(\rmd x) \PP(\BB \in \rmd M)
\ee
(noting that the first moment measure of $\BB$ is $r \wt \Lambda(\rmd x)$, by Prop. \ref{BN0} (ii)).
Write $\PP_{r+i}(\rmd M) = \PP(\BBi\in \rmd M)$ for $i\in\N_0$ and $M \in \mathbb{M}$.
Then, from \eqref{b1a} and \eqref{palm0}, 
\begin{align}\label{b1b}
\PP(\wt J_1 \in \rmd x, \,  \BB \in \rmd M) &=  \frac{x}{T(M)}M(\rmd x)\PP_r(\rmd M) \nonumber \\
&=  \frac{x}{T(M)} r \wt \Lambda(\rmd x) \PP^{(x)}_r(\rmd M).
\end{align}
By \eqref{palmBN}, $\PP_r^{(x)}$ is the distribution of $\delta_x + \xi$ where $\xi$ is distributed as $\BN(r+1, \wt \Lambda)$.
For each $x\in(0,1)$, let $\BBx= \BB - \delta_x$. 
Changing variable to $M_1 = M- \delta_x$ in \eqref{b1b} gives 
\be\label{b1c}
\PP(\wt J_1 \in \rmd x, \,  \BBx\in \rmd M_1)=  
\frac{x}{T(M_1)+x} r \wt \Lambda(\rmd x) \PP_{r+1}(\rmd M_1).
\ee
Then  noting that, jointly,  $ \PP(\wt J_1 \in \rmd x, \,  \BB_1 \in \rmd M_1)= \PP(\wt J_1 \in \rmd x, \,  \BBx\in \rmd M_1)$, we have proved \eqref{b1}.

(ii)\ 
$ ^{(r)}T_1=T(\BB_1)$ is a deterministic transformation of $\BB_1$, so for each $y>0$,
\bean\label{b1Ta}
\PP(\wt J_1 \in \rmd x,\, \BB_1 \in \rmd M, \, T(\BB_1)\le y) 
&=&
\PP(\wt J_1 \in \rmd x,\, \BB_1 \in \rmd M,\, M\in Q_y)\cr
&=&
{\bf 1}_{\{M\in Q_y\}}\PP(\wt J_1 \in \rmd x,\, \BB_1 \in \rmd M),
\eean
where $Q_y:= \{M\in \mathbb{M}:T(M)\le y\}$.
By \eqref{b1} the last expression equals 
\ben
\frac{x}{T(M) + x} r \wt \Lambda(\rmd x ) 
\PP \big( \BBy \in \rmd M, \, T(\BBy) \le y\big),
\een
from which \eqref{b1T} follows. 

(iii)\
Integrating $M$ out of \eqref{b1T} and recalling \eqref{gd} gives \eqref{b1T0} via 
\ben\label{b1Tb}
\PP(\wt J_1 \in \rmd x, \,  \,  T(\BB_1) \in \rmd t  ) 
= \frac{x}{t + x} r \wt \Lambda(\rmd x ) 
\PP\big(T(\BBy) \in \rmd t \big).		
\een
This completes the proof of Proposition \ref{b11}. \halmos

We can now compute the joint density of the size-biased points of $\BN(r, \wt \Lambda)$.
Write the ascending factorial 
as $r^{(n)} = r(r+1)\cdots(r+n-1)$, $n\in\N$.

\begin{proposition}\label{joint_density}
Fix $r, n\in\N$. 
Given $x_i \in (0,1) $, $1\le i\le n$, $x_i \neq x_j$ for $i \neq j$,  and $t > \sum_{i=1}^n x_i$, 
we have (interpreting  $\sum_0^1\equiv 0$)
\bea\label{delta_jd}
&&
\PP\big(\wt J_1 \in \rmd x_1,  \ldots, \wt J_n \in \rmd x_n, \,  {}^{(r)}T \in \rmd t\big)\cr
&& = r^{(n)}\alpha^n 
 \prod_{i=1}^n \frac{x_i^{-\alpha} \rmd x_i}{t - \sum_{j=1}^{i-1}x_j} \PP\Big( {}^{(r+n)}T \in \rmd \Big(t - \sum_{i=1}^n x_i \Big) \Big).
\eea
\end{proposition}

\medskip\noindent {\bf Proof of Proposition \ref{joint_density}:}\
Given $x_1, \ldots, x_n \in (0,1)$, $x_i \neq x_j$ for $i \neq j$, and $M \in \mathbb{M}$, write 
$M_{i+1} = M_{i} - \delta_{x_{i+1}}$, 
with $M_0 = M$ and $i=0,\ldots, n-1$.  
We consider only the first two size-biased picks with $x_1 \neq x_2$. The extension to general $n$ is similar.
Letting $\mathbb{M}_1:= \MM- \delta_{x_1}$, we compute
\begin{align}\label{palm1a}
&\PP\big(\wt J_1 \in \rmd x_1,  \, \wt J_2 \in \rmd x_2, \, T(\BB) \in \rmd t \big) \nonumber \\
=& \int_{M \in \mathbb{M}} \PP\big(\wt J_1 \in \rmd x_1,  \, \wt J_2 \in \rmd x_2, \, \BB \in \rmd M, \,  T(\BB) \in \rmd t \big) \nonumber \\
=& \int_{M_1 \in \mathbb{M}_1} \PP\big(\wt J_1 \in \rmd x_1,  \, \wt J_2 \in \rmd x_2, \, 
\BX \in \rmd M_1, \,  T(\BX) \in \rmd (t-x_1) \big) \nonumber\\
=& \int_{M_1 \in \mathbb{M}_1}   \PP\big( \wt J_1 \in \rmd x_1,  \, \wt J_2 \in \rmd x_2, \, \BB_1 \in \rmd M_1, \, T(\BB_1) \in \rmd (t-x_1) \big).
\end{align}

The probability on the RHS of \eqref{palm1a} can be replaced by
\begin{align}\label{jd_a}
&\,   \PP\big(\wt J_2 \in \rmd x_2 \,\big |  \, \wt J_1 = x_1, \, 
\BB_1 = M_1,  {}^{(r)}T_1 = t-x_1 \big)
 \nonumber \\
 &\hskip2cm 
 \times \PP\big(\wt J_1 \in \rmd x_1, \, 
 \BB_1 \in  \rmd  M_1, {}^{(r)}T_1 \in \rmd (t-x_1) \big) \nonumber \\
=& \, \PP\big(\wt J_2 \in \rmd x_2 \,\big |  \, \BB_1 = M_1,  {}^{(r)}T_1 = t-x_1 \big) 
\nonumber \\
&\hskip2cm 
\times
\frac{x_1}{t}r \wt \Lambda(\rmd x_1) 
\PP\big(\BBy \in \rmd M_1,\, T(\BBy) \in \rmd (t-x_1) \big)\nonumber \\
=& 
\,\frac{x_2}{ t - x_1} M_1(\rmd x_2) \nonumber \\
&\hskip2cm \times 
\frac{x_1}{t}r \wt \Lambda(\rmd x_1) \PP\big(\BBy \in \rmd M_1,\, T(\BBy) \in \rmd (t-x_1) \big).
\end{align}
The first equality in \eqref{jd_a} comes from \eqref{b1T} and the fact that $\wt J_2$ is conditionally independent of $\wt J_1$ given $\BB_1$. In the last equality of \eqref{jd_a}, we used the definition of  size-biased picks, as in \eqref{b1a}.
Using \eqref{jd_a}, the RHS of \eqref{palm1a} equals
\begin{align}\label{jd_d}
 \frac{x_2}{ t - x_1}\frac{x_1}{t}r  \wt \Lambda(\rmd x_1)  
 \int_{M_1 \in \mathbb{M}_1}M_1(\rmd x_2) \PP\big(\BBy \in \rmd M_1,
 T(\BBy) \in \rmd (t-x_1) \big). 
\end{align}

When $\xi$ is a point process with distribution $\PP_{r+1}^{(x)}$, for each $M \in \MM$ and $t > 0$, we abbreviate $ \PP\big( \xi \in \rmd M, \, T(\xi) \in \rmd t \big)$ to  $\PP_{r+1}^{(x)}(\rmd M, \rmd t) $. 
Recalling that ${}^{(r)}T_1 = T(\BB_1)$ is a deterministic function of $\BB_1$, then by \eqref{palm0}, 
\[ M_1(\rmd x_2) \PP_{r+1}\big(\rmd M_1, \rmd (t-x_1) \big)= 
(r+1) \wt \Lambda(\rmd x_2) \PP_{r+1}^{(x_2)}\big(\rmd M_1,  \rmd (t-x_1) \big) .
\]
Substituting this in \eqref{jd_d}, we obtain from \eqref{palm1a} that
\begin{align*}  
&\PP\big(\wt J_1 \in \rmd x_1,  \, \wt J_2 \in \rmd x_2, \, T(\BB) \in \rmd t \big) \nonumber \\
&=  \frac{x_2}{ t - x_1}\frac{x_1}{t}r  \wt \Lambda(\rmd x_1)  
 \int_{M_1 \in \mathbb{M}_1}(r+1) \wt \Lambda(\rmd x_2) \PP_{r+1}^{(x_2)}\big(\rmd M_1,  \rmd (t-x_1) \big) \nonumber \\
&=  \frac{x_2}{ t - x_1}\frac{x_1}{t}r \wt \Lambda(\rmd x_1) (r+1) \wt \Lambda (\rmd x_2)  \cr
& \hskip3cm
\times \int_{M_2 \in \mathbb{M}_2}  \PP \big(\BBz \in \rmd M_2, T(\BBz) \in \rmd (t-x_1-x_2) \big)\, \nonumber \\
&= \, \prod_{i=1}^2 \frac{(r+i-1)x_i\wt \Lambda(\rmd x_i)}{t - \sum_{j=1}^{i-1}x_i}   \PP \big( {}^{(r+2)}T\in \rmd (t-x_1-x_2) \big).
\end{align*}
Here the second equality is obtained by changing variable to $M_2 = M_1 - \delta_{x_2}$  as in \eqref{b1c}, with $\mathbb{M}_2:= \MM_1 - \delta_{x_2}$.
In the last equality, we note that $\PP(\BBz \in \MM_2) =1$, as $\BBz$ has a diffuse mean measure, hence
\[
\PP\big(\BBz(\{x_1\}) > 0\big) = \PP\big(\BBz (\{x_2\}) > 0\big) = 0.
\]

This proves \eqref{delta_jd} when $n = 2$.
By a similar argument, we can show that for each $n \in \NN$,
$x_i\in(0,1)$, $t>\sum_1^nx_i$, 
\bean
&&\PP\big( \wt J_1 \in \rmd x_1, \, \ldots,  \,  \wt J_n \in \rmd x_n,\, T(\BB) \in \rmd t \big) \nonumber \\
&&=  \prod_{i=1}^n \frac{ (r+i-1) x_i \wt \Lambda(\rmd x_i)}{t - \sum_{j=1}^{i-1}x_j} 
\PP\Big( {}^{(r+n)}T \in \rmd \Big(t - \sum_{i=1}^n x_i \Big)\Big),
\eean
and this is the same as \eqref{delta_jd}. \halmos

\medskip
Next we use \eqref{delta_jd} to derive the joint densities of the size-biased quantities in \eqref{sb_T}--\eqref{Un}.
Write $\Theta (x) =  \alpha x^{-\alpha} {\bf 1}_{\{0 <x <1\}}$
and recall ${}^{(r)}T_0\equiv {}^{(r)}T$.

\begin{proposition}\label{jointTT}
Fix $r \in\NN$.

{\rm (i)}\ The joint density of $\left({}^{(r)}T, \, {}^{(r)}T_1, \, {}^{(r)}T_2, \ldots,  {}^{(r)}T_n\right)$ with respect to Lebesgue measure is, for $t_0 >t_1> \cdots > t_{n} > 0$  and $n \in \NN$,
\be\label{jointT}
f(t_0, t_1, \ldots, t_n) = r^{(n)}  g_{r+n}(t_n) \prod_{i=0}^{n-1} \frac{\Theta(t_i - t_{i+1})}{t_i}.
\ee

{\rm (ii)}\
The sequence $\left({}^{(r)}T, \, {}^{(r)}T_1, \, {}^{(r)}T_2, \ldots \right)$ is a (non-homogeneous) Markov Chain with transition density, for $t_n>t_{n+1}>0$ and  $n \in \N_0$,
\be\label{transit}
\PP\big(  {}^{(r)}T_{n+1} \in  \rmd t_{n+1} \, \big| \, {}^{(r)}T_n = t_n   \big)=  (r+n)  \frac{\Theta(t_n-t_{n+1})}{t_n}  \frac{g_{r+n+1}(t_{n+1})}{g_{r+n}(t_n)}  \rmd t_{n+1}.
\ee

{\rm (iii)}\ The joint density of $\left({}^{(r)}T_n, \, {}^{(r)}U_1, \, {}^{(r)}U_2, \ldots,  {}^{(r)}U_n\right)$ is, for $t_n>0$, $0<u_i<1$, $1\le i\le n$,  and $n \in \NN$,
\begin{align}\label{jointUTn}
h(t_n, u_1, \ldots, u_n) 
=&
\frac{r^{(n)}}{K_{n}}  g_{r+n} (t_n)  t_n^{-n\alpha} \nonumber \\
&\times 
 \prod_{i=1}^n \frac{\Gamma(i\alpha+1-\alpha)}{\Gamma(i\alpha)\Gamma(1-\alpha)} u_i^{i\alpha-1} \bar u_i^{-\alpha} {\bf 1}_{\{ t_n <\prod_{j=i}^n u_j/\bar u_i\}}, 
\end{align}
where $\bar u_i = 1 - u_i$, and 
\be\label{Kn}
 K_{n}
= 
\frac{\prod_{i=0}^{n-1} \Gamma(1+i\alpha)}
{\alpha^{n} \Gamma^n(1-\alpha) \prod_{i=1}^n \Gamma(i\alpha)} 
=
\frac {\Gamma(n+1)}{\Gamma^n(1-\alpha) \Gamma(n\alpha+1)}.
\ee
\end{proposition}

\begin{remark}
{\rm 
A routine calculation shows that
\be\label{Kn1}
\frac{1}{\Gamma(n\alpha)}\int_0^\infty t^{n\alpha-1} \Big( 1+ \int_0^\infty (1-e^{-tx})\alpha x^{-\alpha-1} \rmd x\Big)^{-r-n}\rmd t
=\frac{K_n}{r^{(n)}}.
\ee
Notice that the inner integration in \eqref{Kn1} is over $x\in(0,\infty)$, whereas that in \eqref{gdL} is over $x\in(0,1)$, and our $\Theta(x)$ is restricted to $(0,1)$, whereas that of PPY \citeyearpar[Eq. (2.b)]{PPY1992} is not. 
This is a reflection of the truncation induced by eliminating the large points. Still, the 
$K_n$ in \eqref{Kn} and \eqref{Kn1} exactly equals the $K_n$ in Eq. (2.n) of PPY~\citeyearpar{PPY1992}, when the stable scaling constant $c$ in their Eq. (2.i) is set equal to $\Gamma(1-\alpha)$.
In both notations,  $K_n = \EE(S_1^{-n\alpha})$ 
(see also Eq. (30) of PY \citeyearpar{PY1997}).

In general, we have the following relation:
\begin{align}\label{depen}
K_{n} &= r^{(n)}\int_{u_1=0}^1 \cdots \int_{u_n=0}^1 \int_{t_n = 0}^{d(u_1, \ldots, u_n)} 
t_n^{-n\alpha} g_{r+n}(t_n) \rmd t_n \nonumber \\
&\hspace{2in} \times \prod_{i=1}^n f_{B_{i\alpha, 1-\alpha}}(u_i)\rmd u_1 \cdots \rmd u_n, 
\end{align}
where $d(u_1, \ldots, u_n):= \min_{1\le i \le n} \prod_{j=i}^n u_j/\bar u_i$ 
for  $0<u_i<1$, $1\le i\le n$,  $n \in \NN$, and
 $f_{B_{a,b}}$ is the density of a Beta(a,b) distribution as in \eqref{betadef}.
 }

\end{remark}

\medskip\noindent {\bf Proof of  Proposition \ref{jointTT}:}\
(i):\ 
By change of variable in \eqref{delta_jd}, we have
\bean 
 &&
 \PP \big({}^{(r)}T \in \rmd t_0, \, {}^{(r)}T_1 \in \rmd t_1, \, {}^{(r)}T_2 \in \rmd t_2, \ldots,  {}^{(r)}T_n \in \rmd t_n\big) 
 \nonumber \\
&&=\, \PP \big(  {}^{(r)}T \in \rmd t_0, \, \wt J_1 \in \rmd (t_0-t_1),\, \wt J_2 \in \rmd (t_1-t_2),  \ldots, \wt J_n \in \rmd (t_{n-1} - t_n) \,\big) \nonumber \\
&&=\, r^{(n)}  \prod_{i=0}^{n-1} \frac{\Theta(t_i - t_{i+1})}{t_i} g_{r+n}(t_n) \, \rmd t_0 \rmd t_1 \cdots \rmd t_n.
\eean
This proves \eqref{jointT}.
Part (ii) follows immediately from Part (i):
\bean  
&&\, \PP\big({}^{(r)}T_{n+1} \in \rmd t_{n+1} \, \big|\, {}^{(r)}T = t_0, \, {}^{(r)}T_1= t_1, \, {}^{(r)}T_2 = t_2, \ldots,  {}^{(r)}T_n = t_n \big) \nonumber \\
&&=\,  (r+n)  \frac{\Theta(t_n-t_{n+1})}{t_n}  \frac{g_{r+n+1}(t_{n+1})}{g_{r+n}(t_n)}\rmd t_{n+1},
\eean
which does not depend on $t_0,t_1,\ldots, t_{n-1}$.
Thus \eqref{transit} is established.

(iii)\ To show \eqref{jointUTn}, we first consider the case $n=2$.
Note that 
\ben 
h(t_2, u_1,  u_2) = f\Big(\frac{t_2}{u_1 u_2 },  \frac{t_2}{u_2},t_2 \Big) t_2^2u_1^{-2} u_2^{-3};
\een
where $f$ is defined in \eqref{jointT} and $t_2^2u_1^{-2} u_2^{-3}$ is the 
Jacobian 
from the change of variables.  Expanding the expression in \eqref{jointT} with $\Theta(x) = \alpha x^{-\alpha}{\bf 1}_{\{0< x < 1\}}$, we get $h(t_2, u_1, u_2)$ equal to
\begin{align*}  
 &r^{(2)} g_{r+2}(t_2) u_1^{-1} u_2^{-1} \Theta\bigg(\frac{t_2}{u_1u_2} \bar u_1\bigg) \Theta\bigg(\frac{t_2}{u_2} \bar u_2\bigg) \nonumber \\
=\, & r^{(2)} g_{r+2}(t_2) \alpha^2 t_2^{-2\alpha}  \big(u_2^{2\alpha-1}\bar u_2^{-\alpha}\big) \big( u_1^{\alpha-1} \bar u_1^{-\alpha} \big){\bf 1}_{\{ t_2 \bar u_2/ u_2< 1\}}{\bf 1}_{\{ t_2 \bar u_1/ (u_1u_2)< 1\}} \nonumber \\
=\, & \frac{r^{(2)}}{K_{2}} \, g_{r+2}(t_2)  t_2^{-2\alpha} \cdot  \bigg[ \frac{\Gamma(1+\alpha)}{\Gamma(2\alpha) \Gamma(1-\alpha)}u_2^{2\alpha-1}\bar u_2^{-\alpha}\bigg] 
\nonumber \\
&\qquad \qquad \qquad \times
\bigg[ \frac{\Gamma(1)}{\Gamma(\alpha) \Gamma(1-\alpha)}  u_1^{\alpha-1} \bar u_1^{-\alpha} \bigg]{\bf 1}_{\{ t_2 < u_2/\bar u_2\}}{\bf 1}_{\{ t_2 < u_1u_2/\bar u_1\}},
\end{align*}
where 
\[ 
K_{2} =\frac { \prod_{i=0}^1 \Gamma(1+i\alpha)}
{\alpha^{2}\Gamma^2(1-\alpha)\prod_{i=1}^2 \Gamma(i\alpha)}.
\]
This formula can be generalised to $n\ge 2$ similarly, and \eqref{jointUTn} follows.    \halmos

\medskip
To complete this section our final theorem gives formulae for the distributions of the size-biased  sequence constructed from $\PD_\alpha^{(r)}$ as defined in 
in \eqref{sub_alpha}, as well as for the residual fractions defined in \eqref{U1}--\eqref{Un}.

\begin{theorem}\label{stick}
{\rm (i)}\ For each $r\in\N$ let $(V_n^{(r)})_{n\in\N}$ have a $\PD_\alpha^{(r)}$ distribution as defined in
\eqref{sub_alpha}, with corresponding size-biased sequence $(\wt V_n^{(r)})$.
Then for each $n\in\N$ the joint density of 
$\big(\wt V_1^{(r)},  \ldots, \wt  V_n^{(r)},  {}^{(r)}T\big) $ with respect to Lebesgue measure is
\be\label{sss_jd}
p_r(v_1, \ldots, v_n,t) = 
r^{(n)}\alpha^n \times
\prod_{i=1}^n \frac{ v_i^{-\alpha} }{\bar  v_{i-1}} {\bf 1}_{\{v_i < 1/t \}}\times t^{-n\alpha} g_{r+n}(t\bar v_n),
\ee
where $t>0$, $0<v_i< 1$ are such that $ \sum_{i=1}^n v_i< 1$,  
  $\bar v_0\equiv 1$, and, for each $i\ge 1$, $\bar v_i= 1 -v_1 - \cdots- v_i$.

{\rm (ii)}\ The joint distribution  of ${}^{(r)} T_n$ and
${}^{(r)}U_1, \, {}^{(r)}U_2, \ldots,  {}^{(r)}U_n$ 
can be written as 
\be\label{h6}
\big({}^{(r)}T_n, \, {}^{(r)}U_1, \, {}^{(r)}U_2, \ldots,  {}^{(r)}U_n\big)
\eqd \big( Y_{d(U_1, \ldots, U_n)},  U_1,  U_2, \ldots,  U_n\big)å,
\ee
where the  $(U_i)$ are independent Beta$(i\alpha, 1-\alpha)$ rvs,
independent of $\BB$,   the function 
$d(u_1, \ldots, u_n):= \min_{1\le i \le n} \prod_{j=i}^n u_j/\bar u_i$ 
and, for each $c > 0$, 
$Y_c \eqd ({}^{(r+n)}T)^{-n\alpha}{\bf 1}_{\{{}^{(r+n)}T < c\}}$.
\end{theorem}

\medskip\noindent {\bf Proof of  Theorem \ref{stick}:}\
(i)\  
Identify the size-biased $\wt V_i^{(r)}$ with the points 
$\wt J_r(i)$ in \eqref{Jd} normalised by their sum
$T(\BB)$. 
Then change variable in \eqref{delta_jd} to $v_i = x_i/t$ and substitute for $\wt \Lambda$ to get 
\begin{align*}
&\PP\big( \wt V^{(r)}_1 \in \rmd v_1, \, \ldots,  \,  \wt V^{(r)}_n \in \rmd v_n,\, T(\BB) \in \rmd t \big)  \nonumber \\
=&
\PP\big( \wt J_r(1) \in t\rmd v_1,\, \ldots,\, \wt J_r(n)\in t\rmd v_n,\,  {}^{(r)}T  \in \rmd t \big)  \nonumber \\
=&
  r^{(n)}\alpha^{n} \prod_{i=1}^n \frac{t\cdot t^{-\alpha}v_i^{-\alpha}
  \rmd v_i }{t - \sum_{j=1}^{i-1} t v_j}  \times
\PP\Big( {}^{(r+n)}T \in \rmd \Big(t - \sum_{i=1}^n x_i \Big)\Big)
\quad {\rm (by}\ \eqref{delta_jd})
\nonumber \\
=&
 r^{(n)}\alpha^{n}  t^{-n\alpha} \prod_{i=1}^n \frac{ v_i^{-\alpha} \rmd v_i }{\bar v_{i-1}}
\PP\big( {}^{(r+n)}T \in \rmd t \bar v_n \big),
\end{align*}   
Recall that $g_{r+n}(\cdot)$ is the density of $^{(r+n)}T$ (see \eqref{gd}),  to complete the proof of \eqref{sss_jd}.

(ii)\ 
Representation 
\eqref{h6} is immediate from Proposition \ref{jointTT} and \eqref{depen},
with the vector on the RHS of \eqref{h6} having the structure specified.
\halmos

\begin{remark}\label{rem3d} 
{\rm 
(i)\
There are some quite involved manipulations in obtaining the above formulae. As a check on the calculations, in the Appendix of  \href{https://arxiv.org/abs/1611.09980}{arxiv:1611.09980}, we give direct verifications that \eqref{jointT} and  \eqref{sss_jd} are probability densities  (integrate to 1). 

(ii)\ 
For a stick-breaking representation, solve \eqref{U1} and \eqref{Un} to get
\be\label{wVsb}
\wt V_n^{(r)} = (1- {}^{(r)}U_n) \prod_{i=1}^{n-1} {}^{(r)}U_i.
\ee
The joint distribution of $({}^{(r)}U_i)_{1\le i\le n}$ can be computed from \eqref{h6}, in which we note that $U_1,  U_2, \ldots,  U_n$ are individually independent 
but dependence overall is introduced via the connection with the $Y$ term.
In this respect the result is different from the $\PD_\alpha$  situation, as we would expect, but the distribution of $\wt V_n^{(r)}$ as given by \eqref{wVsb} is sufficiently  explicit to enable computations or simulations. 

\rm (iii)\ \eqref{sss_jd} generalises the corresponding version for $\PD_\alpha={\rm PD}(\alpha,0)$ in PY~\citeyearpar[Prop. 47]{PY1997}.

(iv)\
When we sample from a Poisson process,  the various quantities in Proposition \ref{jointTT} are computed in PPY \citeyearpar[Theorem 2.1]{PPY1992}. 

(v)\ Although motivated by the idea of trimming an integer number $r$ of large jumps, our formulae once derived are valid for $r>0$, and available for modelling purposes in this generality.
}
\end{remark}
%
To conclude this section  we expand briefly on the differences as well as the similarities between the $\PD_\alpha$ and
 $\PD_\alpha^{(r)}$  approaches.
 In both cases, start with a stable($\alpha$) subordinator $S$
with ranked jumps $\Delta S_1^{(1)}\ge \Delta S_1^{(2)}\ge\cdots$.
The sequence $(\Delta S_1^{(i)}/S_1)_{i\ge 1}$ then has a $\PD_\alpha$ distribution. 
We can think of these as the points from a Poisson point process
with intensity measure $\Lambda(\rmd x)=\alpha x^{-\alpha-1}\rmd x$, normalised by their sum.
For  $\PD_\alpha^{(r)}$, the analogous process is the negative binomial point process   $\BN(r, \wt \Lambda)$ 
formed from  ratios of jumps rather than from the jumps themselves, i.e., 
\ben 
 \BB = \sum_{i \ge 1} \delta_{J_r(i)}, \quad \text{with} \quad 
 J_r(i) = \frac{\Delta S_1^{(r+ i)}}{\Delta S_1^{(r)}},\ i\in\N.
\een
The normalised jumps  on which a size-biased version is based are
\be\label{nj}
 \frac{J_r(i)}{\sum_{\ell\ge 1}J_r(\ell)}
=
 \frac{\Delta S_1^{(r+ i)}}{\Delta S_1^{(r)}}
 \bigg/\frac{^{(r)}S_1}{\Delta S_1^{(r)}}
 = \frac{\Delta S_1^{(r+i)}}{ ^{(r)}S_1},\ i\in\N,
\ee
and the sequence formed from  these has a $\PD_\alpha^{(r)}$ distribution, as we define it,   on the infinite simplex. 

We may set $r=0$ in \eqref{sub_alpha}  to have the distribution   of  $(V_n^{(r)})_{n\in\N}$, that is,  $\PD_\alpha^{(r)}$, reduce to that of  $(V_n)_{n\in\N}$, that is, $\PD_\alpha$. 
But we cannot take $r=0$ in \eqref{nj} 
with the idea that the size-biased distribution associated with 
$\PD_\alpha^{(r)}$ might then reduce to the one associated with $\PD_\alpha$.
Our analysis proceeds via the process $\BB$,  which is not defined for $r=0$ (its points $J_r(i)$ are not defined for $r=0$).  
Setting $r=0$ in formulae such as \eqref{jointT}, \eqref{jointUTn}, \eqref{sss_jd}, etc., which result from an analysis of $\BB$,
is not permissible.

%

\section{Discussion and Applications}\label{dis}
We mention two applications papers which vividly illustrate the possibilities for useful and revealing application of our results.
 \cite{sos15} shows ``capital distribution curves" (CDCs) for
over 20 countries listed on the NASDAQ stock exchange. 
 The CDC  is a log plot of normalized stock capitalizations ranked in descending order, against their log-ranks. The curves display remarkable stability over periods of time and are very well fitted by a $\PD_\alpha$ distribution over much of their range.
But a glance at  Sosnovskiy's Figure 3, for example,  shows that an even better fit would result from discarding a small number of the largest stocks -- 
suggesting a $\PD_\alpha^{(r)}$ distribution.
We might indeed expect that a small  number of very large stocks would show aberrant behaviour, compared to the majority.

A very similar situation occurs with the Zipf plots (log frequencies of words in the Penn Wall St. journal, against their log-ranks), in the paper of  \cite{ggj11}.
Half a dozen or so of the most frequent words appear as outliers, while the rest conform closely to a $\PD_\alpha$ fit (see their Figure 4).
 
 In general, we can expect that our generalised $\PD^{(r)}_\alpha$ distribution could be used to robustify analyses and reveal interesting features in this kind of data.

\medskip\noindent {\bf Acknowledgements.}\
We are grateful for some very helpful feedback from 
Peter  Kevei and David Mason.
 
\section{Appendix: (\ref{jointT}) and (\ref{sss_jd})  integrate to 1}\label{sA}\
As a check on the calculations, we give here a direct verification that \eqref{jointT} and  \eqref{sss_jd} integrate to 1  in the case $n=1$. 
An extension to larger $n$ is straightforward.

Eq. \eqref{jointT} gives, for $n=1$, $r\in\N$, $t_0>t_1>0$,
\bean 
f(t_0, t_1) = r g_{r+1}(t_1)\frac{\Theta(t_0-t_{1})}{t_0}
=r\alpha g_{r+1}(t_1)\frac{(t_0-t_{1})^{-\alpha}{\bf 1}_{\{t_0-t_1<1\}}}{t_0}
\eean
Notice that
\ben
\int_{t_0=t_1}^{1+t_1} \frac{(t_0-t_{1})^{-\alpha}}{t_0}\rmd t_0
=\int_0^1 \frac{t^{-\alpha}}{t+t_1}\rmd t
\een
so
\bea\label{jT}
\int_{t_1=0}^\infty \int_{t_0=t_1}^{1+t_1} f(t_0, t_1) \rmd t_1\rmd t_0
&=&
r\alpha\int_{t_1=0}^\infty  g_{r+1}(t_1)
\int_{t=0}^1 \frac{t^{-\alpha}}{t+t_1}\rmd t\cr
&=&
r\alpha\int_{t=0}^1 t^{-\alpha} \rmd t \int_{t_1=0}^\infty 
\frac{g_{r+1}(t_1)}{t+t_1}\rmd t_1.
\eea
Introduce an integral over $\lambda$ and then substitute from \eqref{gd} to write the last expression as
\bean
&&
r\alpha\int_{t=0}^1 t^{-\alpha} \rmd t \int_{t_1=0}^\infty \int_{\lambda=0}^\infty e^{-\lambda(t+t_1)} \rmd \lambda\,
g_{r+1}(t_1)\rmd t_1\cr
&&=
r\alpha\int_{t=0}^1 t^{-\alpha} \rmd t \int_{\lambda=0}^\infty e^{-\lambda t}
\int_{t_1=0}^\infty e^{-\lambda t_1}
g_{r+1}(t_1)\rmd t_1\,  \rmd \lambda\cr
&&=
r\alpha\int_{t=0}^1 t^{-\alpha} \rmd t \int_{\lambda=0}^\infty 
\frac{e^{-\lambda t} \rmd \lambda}
{\Big(1 + \int_{0}^{1} (1- e^{-\lambda x}) \wt \Lambda(\rmd x)  \Big)^{r+1}}.
\eean
The last equation used \eqref{gdL}. 
The final integral can be evaluated as
\bean
r\int_{\lambda=0}^\infty 
\frac{\alpha\int_{t=0}^1 t^{-\alpha}  e^{-\lambda t} \rmd t\, \rmd \lambda}
{\Big(1 + \int_{0}^{1} (1- e^{-\lambda x}) \wt \Lambda(\rmd x)  \Big)^{r+1}}
&=&
r
\frac{\Big(1 + \int_{0}^{1} (1- e^{-\lambda x}) \wt \Lambda(\rmd x)  \Big)^{-r}}{-r}\Bigg|_0^\infty \cr
&=&
1.
\eean

Next we give a direct verification that \eqref{sss_jd} integrates to 1  in the case $n=1$. 
\bean\label{jT1}
&&
r\alpha\int_{v_1=0}^1 v_1^{-\alpha}
\int_{t=0}^{1/v_1} t^{-\alpha}  g_{r+1}(t(1-v_1))\rmd t\cr
&&=
r\alpha\int_{v_1=0}^1 v_1^{-\alpha} (1-v_1)^{\alpha-1}
\int_{t=0}^{(1-v_1)/v_1} t^{-\alpha}  g_{r+1}(t)\rmd t\cr
&&=
r\alpha\int_{y=0}^\infty \frac{y^{\alpha-1}}{1+y}\rmd y
\int_{t=0}^{y} t^{-\alpha}  g_{r+1}(t)\rmd t
\quad (y=(1-v_1)/v_1)
\cr
&&=
r\alpha\int_{t=0}^{\infty} t^{-\alpha}  g_{r+1}(t)\rmd t
\int_{y=t}^\infty \frac{y^{\alpha-1}}{1+y}\rmd y\cr
&&=
r\alpha\int_{t=0}^{\infty} g_{r+1}(t)\rmd t
\int_{y=1}^\infty \frac{y^{\alpha-1}}{1+ty}\rmd y
\quad (y=y/t)\cr
&&=
r\alpha\int_{y=1}^\infty y^{\alpha-1}\rmd y
\int_{t=0}^{\infty} \frac{1}{1+ty} g_{r+1}(t)\rmd t
\cr
&=&
r\alpha\int_{x=0}^1 x^{-\alpha} \rmd x
\int_{t=0}^\infty \frac{g_{r+1}(t)}{x+t}\rmd t
\quad (x=1/y).
\eean
The last is the RHS of \eqref{jT} which equals 1.

\bibliography{Library_Levy_Mar2017}
\bibliographystyle{newapa}
\end{document}